\documentstyle[amscd,12pt,psamsfonts,leqno,amssymb]{amsart}
\input xypic

\newtheorem{defn}{Definition}
\newtheorem{thm}[defn]{Theorem}

\newtheorem{cor}[defn]{Corollary}
\newtheorem{lem}[defn]{Lemma}
\newtheorem{prop}[defn]{Proposition}

\theoremstyle{remark}
\newtheorem{rem}{Remark}
\theoremstyle{remark}

\numberwithin{equation}{section}
\numberwithin{defn}{section}

\renewcommand\sp{\operatorname{Spec}}

\newcommand\grv{{\operatorname{Gr}}(V)}
\newcommand\gr{\operatorname{Gr}}

\newcommand\gl{\operatorname{Gl}}
\newcommand\lgl{\operatorname{LGl}}
\renewcommand\hom{\operatorname{Hom}}

\newcommand\Det{\operatorname{Det}}

\newcommand\res{\operatorname{Res}}

\newcommand\limil[1]{\underset{#1}\varinjlim\,}

\newcommand\limpl[1]{\underset{#1}\varprojlim\,}
\newcommand\rad{\operatorname{Rad}}
\newcommand\aut{\operatorname{Aut}}
\renewcommand\hom{\operatorname{Hom}}

\newcommand\tr{\operatorname{Tr}}
\newcommand\Detd{\operatorname{Det}^\ast}

\newcommand\id{\operatorname{Id}}

\newcommand\der{\operatorname{Der}}

\newcommand\pic{\operatorname{Pic}}
\newcommand\lie{\operatorname{Lie}}

\newcommand\kr{\operatorname{Kr}}

\renewcommand\o{{\mathcal O}}

\renewcommand\L{{\mathcal L}}
\renewcommand\H{{\mathcal H}}
\newcommand\M{{\mathcal M}}

\newcommand\Z{{\mathbb Z}}

\newcommand\C{{\mathbb C}}

\newcommand\T{{\mathbb T}}

\newcommand\w{\widehat}
\newcommand\wtilde{\widetilde}

\newcommand\iso{@>{\sim}>>}
\newcommand\beq{
                     \setcounter{equation}{\value{defn}}\addtocounter{defn}1
                     \begin{equation}}

\begin{document}

%%%%%%%%%%%%%%%%%%%%%%%%%%%%%%%%%%%%%%%%%%%%%%%%%%%%%%%%%%%%%
\title{Virasoro Groups and Hurwitz Schemes}
\author[J. M. Mu\~noz Porras and F. J. Plaza Mart\'{\i}n]{J. M. Mu\~noz Porras \\  F. J. Plaza Mart\'{\i}n}
\address{Departamento de Matem\'aticas, Universidad de Salamanca,  Plaza
   de la Merced 1-4
   \\
   37008 Salamanca. Spain.
   \\
    Tel: +34 923294460. Fax: +34 923294583}

\thanks{
    (Secondary).
    \indent This work is partially supported by the research contracts
    BFM2000-1327  of DGI and  SA064/01 of JCyL. The second author is
    also supported by MCYT ``Ram\'on y Cajal'' program.}
\email{jmp@@usal.es} \email{fplaza@@usal.es}

%
%
%\email{esteban@@usal.es} \email{jmp@@usal.es}
%\email{fplaza@@usal.es}

\begin{abstract}
In this paper we study the Hurwitz scheme in terms of the Sato
Grassmannian and the algebro-geometric theory of solitons. We will
give a characterization, its equations and a show that there is a
group of Virasoro type which uniformizes it.
\end{abstract}

\date\today

\maketitle

%%%%%%%%%%%%%%%%%%%%%%%%%%%%%%%%%%%%%%%%%%%%%%%%%%%%%%%%%%%%%

%\setcounter{tocdepth}1
%\tableofcontents

%%%%%%%%%%%%%%%%%%%%%%%%%%%%%%%%%%%%%%%%%%%%%%%%%%%%%%%%%%%%%
\section{Introduction}

Nowadays, the theory of Hurwitz spaces, which parametrize covers
of curves, has shown to be relevant. It has been applied to the
study of the moduli space of curves (\cite{Fulton}) but it is also
important in enumerative geometry (\cite{OP}).

In those problems one focusses on covers of curves $Y\to X$ with
prescribed ramification data. In this paper, since we use the Sato
Grassmannian, these data will be decorated with formal
trivializations. The corresponding functor will be called the
Hurwitz functor and it will be studied in~\S3. The Krichever map
embeds it into the Sato Grassmannian and its image will be
characterized. In particular, the Hurwitz functor is shown to be
representable and its equations will be given.

Let us be more precise. Let us fix a natural number $n$ and a set
of $r$ partitions of it $E=\{\bar e_1,\dots,\bar e_r\}$, where
$\bar e_i=\{e_1^{(i)},\dots, e_{k_i}^{(i)}\}$. Let us denote by
$V$ the $\C$-algebra $\C(\!(z_1)\!)\times \dots\times
\C(\!(z_r)\!)$ and by $W$ the $V$-algebra
    {\small
    $$
    W\,=\,
    \C(\!(z_1^{1/{e_1^{(1)}}})\!)\times \dots \times
    \C(\!(z_1^{1/{e_{k_1}^{(1)}}})\!)
    \times \dots \times
    \C(\!(z_r^{1/{e_1^{(r)}}})\!)\times \dots \times
    \C(\!(z_r^{1/{e_{k_r}^{(r)}}})\!)
    $$}
Then, the Hurwitz functor parametrizes data $(Y,X,\pi,\bar x,\bar
y, t_{\bar x},t_{\bar y})$ where $\pi:Y\to X$ is a cover of
curves, $\bar x=x_1+\dots+x_r\subset X$ and $\sum_{i,j}
y_j^{(i)}\subset Y$ are divisors such that $\pi^{-1}(x_i)=\sum_j
e_j^{(i)} y_j^{(i)}$ and $t_{\bar x}$ and $t_{\bar y}$ are formal
trivializations along $\bar x$ and $\bar y$ respectively (in
particular, they induce isomorphisms $t_{\bar x}:(\w\o_{X,\bar
x})_{(0)}\overset{\sim}\to V$ and  $t_{\bar y}:(\w\o_{Y,\bar
y})_{(0)}\overset{\sim}\to W$). Let $\H^\infty_E[\bar g,g]$ denote
the subscheme of the Grassmannian of $W$ representing the Hurwitz
functor.

In our previous paper \cite{automorphism} we proved that the
Virasoro group scheme, which was defined as the group representing
the functor of automorphisms of $\C(\!(z)\!)$, uniformizes the
moduli space of curves (see Theorem~4.11 of \cite{automorphism}
for the precise statement).

In \S5, we introduce the group $G_V^W$ as a certain subgroup of
the group of automorphisms of $W$ that induce an automorphism on
$V$. We show that this formal group scheme acts canonically on
$\gr(W)$ leaving $\H^\infty_E[\bar g,g]$ stable. Further, we prove
that this group ``uniformizes'' $\H^\infty_E[\bar g,g]$, more
precisely, we prove that the previous action is locally transitive
(Theorem~\ref{thm:GonHloctran}). The proof of these fact requires
the explicit computation of the tangent space to the Hurwitz
space.

In the last section, we add to the previous data a pair
$(L,\phi_{\bar y})$ consisting of a line bundle on $Y$ together
with a formal trivialization of $L$ along $\bar y$. This functor
is representable by a subscheme of the Grassmannian of $W$ which
will be denoted by $\pic^\infty_E[\bar g,g]$ (see
Definition~\ref{defn:pic}). Let $\Gamma_W$ denote the connected
component of $1$ in the scheme representing the functor of
invertible elements of $W$ (see~\S2). Since $G_V^W$ acts
canonically on $\Gamma_W$ we consider the semidirect product
$G_V^W\ltimes \Gamma_W$. Then we show that the group $G_V^W\ltimes
\Gamma_W$ acts on $\pic^\infty_E[\bar g,g]$ and that this actions
is locally transitive (Theorem~\ref{thm:GGammaactPic}).

%%%%%%%%%%%%%%%%%%%%%%%%%%%%%%%%%%%%%%%%%%%%%%%%%%%%%%%%%%%%%
\section{Vector Grassmannians and generalized $\tau$-functions}\label{sect:1}

We assume the base field to be the field of complex numbers.
However, all our results hold for an algebraically closed field of
characteristic zero.

This section recalls and generalizes some results proved in
\S\S2--3 of~\cite{Hurwitz}.

%%%%%%%%%%%%%%
% \subsection*{Formal Geometry}
%%%%%%%%%%%%%%

Let $V$ be the trivial $\C$-algebra
$\C(\!(z_1)\!)\times\overset{r}\dots \times \C(\!(z_r)\!)$. Let us
fix an integer $n>0$ and a set $r$ partitions of it:
    $$
    E\,=\,\{\bar e_1,\dots,\bar e_r\}
    $$
that is, $\bar e_i=\{e_1^{(i)},\dots, e_{k_i}^{(i)}\}$ and
$n=e_1^{(i)}+ \dots+ e_{k_i}^{(i)}$. Finally, let us denote $\bar
r=k_1+\dots+ k_r$.

Associated to this data we consider the following $V$-algebra
    $$
    W^i\,:=\, \C(\!(z_i^{1/{e_1^{(i)}}})\!)\times \dots \times
    \C(\!(z_i^{1/{e_{k_i}^{(i)}}})\!)
    $$
and
    $$
    W\,:=\, W^1\times \dots\times W^r
    $$

Let $V_+$, $W^i_+$ and $W_+$ denote the subalgebras corresponding
to the power series, that is
    $$
    \begin{aligned}
    & V_+\,:=\, \C[\![z_1]\!]\times\overset{r}\dots\times \C[\![z_r]\!] \\
    & W^i_+\,:=\, \C[\![z_i^{1/{e_1^{(i)}}}]\!]\times \dots \times
    \C[\![z_i^{1/{e_{k_i}^{(i)}}}]\!]
    \\
    & W_+\,:=\, W^1_+\times \dots\times W^r_+
    \end{aligned}
    $$

Analogously to the section 2 of~\cite{Hurwitz}, one introduces
formal group schemes whose rational points are the groups of
invertible elements of $V$ and $W$. Let $\Gamma_V$ and $\Gamma_W$
denote the connected components of the origin of these groups.

Let $R$ be a $\C$-algebra. Recall that the set of $R$-valued
points of $\Gamma_V$ is the set of $r$-tuples
$(\gamma_1,\dots,\gamma_r)\in V\w\otimes_{\C} R$ with
$\gamma_i=\sum_{j}a_j^i z_i^{j}$ where $a_j^i\in\rad(R)$ for $j<0$
and $a_0^i$ are invertible. The set of $R$-valued points of
$\Gamma_W$ are described similarly. Indeed, it is the set of
$r$-tuples $(\bar \gamma_1,\dots,\bar \gamma_r)\in W\w\otimes_{\C}
R$ with
    $$
    \bar \gamma_i=
    \big(\sum_{j}a_j^{(i,1)} z_i^{j/e^{(i)}_1}, \dots,
    \sum_{j}a_j^{(i,k_i)} z_i^{j/e^{(i)}_{k_i}}\big)\,\in \,
    W^i\w\otimes_{\C} R
    $$
where $a_j^{(i,h)}\in\rad(R)$ for $j<0$ and $a_0^{(i,h)}$ are
invertible.

%%%%%%%%%%%%%%
%\subsection*{Infinite Grassmannians}
%%%%%%%%%%%%%%

Let us recall briefly some facts on the theory of infinite
Grassmannian. For the sake of simplicity, we will restrict now the
pair  $(W,W_+)$ although all facts hold for $(V,V_+)$ as well.  It
is well known that there is a $\C$-scheme $\gr(W)$, which is
called infinite Grassmannian, whose set of rational points is
   {\small $$\left\{\begin{gathered}
   \text{subspaces $U\subset W$ such that }U\to W/W_+
   \\
   \text{has finite dimensional kernel and cokernel}
   \end{gathered}
   \right\} \ . $$}
This scheme is  equipped with the determinant bundle, $\Det_W$,
which is the determinant of the complex of $\o_{\gr(W)}$-modules
$$
\L\,\longrightarrow\, W/W_+\,\hat\otimes_{\C}\, \o_{\gr(W)} \, ,
$$
where $\L$ is the universal submodule of $\gr(W)$ and the morphism
is the natural projection. The connected components of the
Grassmannian are indexed by the Euler--Poincar\'e characteristic
of the complex. The connected component of index $m$ will be
denoted by $\gr^m(W)$.

The group $\Gamma_W$ acts by homotheties on $W$, and this action
gives rise to a natural action on $\gr(W)$
    $$\Gamma_W\times\gr(W)\,\longrightarrow\,\gr(W) \ .$$
Furthermore, this action preserves the determinant bundle.

\begin{rem}
The linear group $\gl(W)$ acts on  $\gr(W)$ preserving its
determinant bundle (see~\cite{automorphism}~\S2.2). This fact
implies that there is a natural central extension
    $$
    0 \to {\mathbb G}_m \to \wtilde\gl(W)\to \gl(W)\to 0
    $$
In fact, one has such an extension for every subgroup of these
linear groups.

The loop group of the pair $(V,W)$ is defined as the group scheme
$\lgl(W/V)$ representing the subfunctor of groups
$\underline\lgl(W/V)\subset \gl(W)$ defined by
    $$
    \underline\lgl(W/V)(S)\,:=\,
    \left\{
    \begin{gathered}
    \text{automorphisms of }W\w\otimes_{\C} H^0(S,\o_S) \\
    \text{as $V\w\otimes_{\C} H^0(S,\o_S)$-module}
    \end{gathered}
    \right\}
    $$
Since the functors of points of $\lgl(W/V)$ and $\Gamma_W$ are
subfunctors of $\gl(W)$, one obtains central extensions of group
schemes
    $$
    \begin{gathered}
    0\to {\mathbb G}_m\to \wtilde\lgl(W/V)\to \lgl(W/V) \to 0
    \\
    0\to {\mathbb G}_m\to \wtilde\Gamma_W\to \Gamma_W \to 0
    \end{gathered}
    $$
\end{rem}

These facts allow us to introduce $\tau$-functions and
Baker-Akhiezer functions of points of $\gr(W)$. Let us recall the
definition and some properties of these functions
(\cite{Hurwitz},~\S3).

The determinant of the morphism $\L\to W/W_+\hat\otimes_{\C}
\o_{\gr(W)}$ gives rise to a canonical global section
    $$
    \Omega_+\,\in\, H^0(\gr^0(W),\Detd_W) \ .
    $$
In order to extend this section to $\gr(W)$ (in a non-trivial
way), we fix elements $\{v_m\vert m\in\Z\}$ such that: i) the
multiplication by $v_m$ shifts the index by $m$; and, ii)
$v_m\cdot v_{\bar r-r\cdot n-m}=1$. We define
$\Omega_+(U):=\Omega_+(v_m^{-1}U)$ for $U\in\gr^m(W)$.

Now, the $\tau$-function and BA functions will be introduced
following~\cite{Hurwitz}. Recall that
    {\small
    $$
    W\,=\,
    \C(\!(z_1^{1/{e_1^{(1)}}})\!)\times \dots \times
    \C(\!(z_1^{1/{e_{k_1}^{(1)}}})\!)
    \times \dots \times
    \C(\!(z_r^{1/{e_1^{(r)}}})\!)\times \dots \times
    \C(\!(z_r^{1/{e_{k_r}^{(r)}}})\!)
    $$}
and that $\Gamma_W$ parametrices a certain subgroup of invertible
elements of $W$.  Let $t$ be the set of variables
$(t^{(1,1)},\dots, t^{(1,k_1)},\dots,t^{(r,1)},\dots,
t^{(r,k_r)})$ where $t^{(a,b)}=(t_1^{(a,b)}, t_2^{(a,b)}, \dots)$.
Consider the element of $\Gamma_W$ given by
    $$
    g\,=\,
    \big(1+\sum_{j<0}t^{(1,1)}_j z_1^{j/e^{(1)}_1}, \dots,
    1+\sum_{j<0} t^{(r,k_r)}_j z_r^{j/e^{(r)}_{k_r}}\big)\,\in \,
    \Gamma_W
    $$
Then, the $\tau$-function of $U$, $\tau_U(t)$, is defined by
    $$ \tau_U(t)
    \,:=\, \frac{\Omega_+(gU)}{g \delta_U}
    $$
where $\delta_U$ being a non-zero element in the fibre of
$\Detd_W$ over $U$.

Let $z_\centerdot$ denote
$(z_1^{1/{e_1^{(1)}}},\dots,z_1^{1/{e_{k_1}^{(1)}}},\dots,
z_r^{1/{e_1^{(r)}}}, \dots,  z_r^{1/{e_{k_r}^{(r)}}})\in W$.

Let $(a,b)$ be a pair of natural numbers such that $a\in
\{1,\dots,r\}$ and $b\in\{1,\dots, k_a\}$. Let $(c,d)$ be  another
pair satisfying the corresponding constrains.

The $(a,b)$-th Baker-Akhiezer function of a point $U\in\gr(W)$ is
the $W$-valued function whose $(c,d)$-th entry is given by
    {\small
    $$
    \psi_{a,b,U}^{(c,d)}({z_\centerdot},t)\,:=\,
    \exp\Big(-\sum_{i\geq 1}
    \frac{t_i^{(c,d)}}{z_c^{i/ e^{(c)}_d}}\Big)
    \frac{\tau_{U_{a,b}^{c,d}}(t+[z_c^{1/e^{(c)}_d}])}{\tau_U(t)}
    $$}
where
    \begin{itemize}
    \item $[z_v]:=(z_v,\frac{z_v^2}2,\frac{z_v^3}3,\dots)$,
    \item $t+[z_c^{1/e^{(c)}_d}]:=(t^{(1,1)},\dots, t^{(1,k_1)},\dots,
    t^{(c,d)}+[z_c^{1/e^{(c)}_d}] ,\dots ,t^{(r,1)},\dots,
    t^{(r,k_r)})$,
    \item and  $U_{a,b}^{c,d}:=
    (1,\dots,z_a^{1/e^{(a)}_b},\dots,(z_c^{1/e^{(c)}_d})^{-1},\dots,1)\cdot U$.
    \end{itemize}

The main property of these Baker-Akhiezer functions is that they
can be understood as generating functions for $U$ as a subspace of
$W$, as we recall next.

\begin{thm}[\cite{Hurwitz}]\label{thm:BAgeneragrmV}
Let $U\in \gr^m(W)$. Then
    $$
    \begin{aligned}
    \psi_{a,b,U}({z_\centerdot},t)
    \,=\, &
    v_m^{-1} \cdot (1,\dots,z_a^{1/e^{(a)}_b},\dots,1)
    \cdot
    \\ & \qquad \cdot \sum_{i>0}
    \left(\psi_{a,b,U}^{(1,1)}(z_1^{1/e^{(1)}_1}), \dots,
    \psi_{a,b,U}^{(r,k_r)}(z_r^{1/e^{(r)}_{k_r}})\right)p_{a,b,i,U}(t)
    \end{aligned}
    $$
where
    $$
    \big\{
    \big(\psi_{a,b,U}^{(1,1)}(z_1^{1/e^{(1)}_1}), \dots,
    \psi_{a,b,U}^{(r,k_r)}(z_r^{1/e^{(r)}_{k_r}})\big)
     \,\vert\, a\in \{1,\dots,r\}\; , \; b\in\{1,\dots, k_a\} \big\}
    $$
is a basis of $U$ and $p_{a,b,i,U}(t)$ are functions in $t$.
\end{thm}

%%%%%%%%%%%%%%
% \subsection*{Bilinear Identities}
%%%%%%%%%%%%%%

Now, we will follow~\cite{Hurwitz} to prove generalized Bilinear
identities as well as to generalize the associated hierarchy,
which will be called $E$-KP hierarchy.

Recall that $W=W^1\times\dots\times W^r$. Being $W^i$ a
$\C(\!(z_i)\!)$-algebra, one considers pairing given by
    $$
    \begin{aligned}
    W^i\times W^i & \,\longrightarrow\, \C(\!(z_i)\!)
    \\
    (w_1^{(i)}, w_2^{(i)}) &\,\longmapsto \tr^i(w_1^{(i)}w_2^{(i)})
    \end{aligned}
    $$
where $\tr^i:W^i\to \C(\!(z_i)\!)$ is the trace map.

These give rise to a pairing
    $$
    \begin{aligned}
    T_2\colon & W\times W\,\longrightarrow\, \C \\
    & (w_1,w_2) \,\longmapsto \sum_{i=1}^r \res_{z_i=0}
    \tr(w_1^{(i)}w_2^{(i)})dz_i
    \end{aligned}
    $$
where $w_j=(w_j^{(1)},\dots, w_j^{(r)})$ with respect to the
decomposition $W=W^1\times \dots \times W^r$.

From the separability of $W^i$ as $\C(\!(z_i)\!)$-algebra, it
follows that $T_2$ is a non-degenerate bilinear pairing.
Furthermore, it induces an involution of the Grassmannian
    $$
    \begin{aligned}
    \gr(W)\,&\longrightarrow\,\gr(W) \\
    U \,&\longmapsto \, U^\perp
    \end{aligned}
    $$
where $U^\perp$ is the orthogonal of $U$ w.r.t. $T_2$.  This
involution sends the connected component of index $m$ to that of
index $\bar r- r\cdot n- m$.

Finally, the $(a,b)$-th adjoint Baker-Akhiezer functions of $U$
are defined by
    $$
    \psi_{a,b,U}^*({z_\centerdot},t)\,:=\, \psi_{a,b,U^\perp}({z_\centerdot},-t)
    $$

Now, we state the corresponding generalizations whose proofs are
similar to those given in~\cite{Hurwitz}.

\begin{thm}[Bilinear Identity]\label{thm:bilinearGrW}
Let $U,U'\in\gr^m(W)$ be two rational points lying on the same
connected component. Then, $U=U'$ if and only if the following
condition holds
    $$
    T_2\Big(\frac1{z_{\cdot}}\psi_U(z_{\cdot},t),
    \frac1{z_\cdot}\psi_{U'}^*(z_{\cdot},t')\big)\,=\,0
    $$
\end{thm}

%%%%%%%%%%%%%%%%%%%%%%%%%%%%%%%%%%%%%%%%%%%%%%%%%%%%%%%%%%%%%%%%%%%%%5
\section{Hurwitz schemes}
%%%%%%%%%%%%%%%%%%%%%%%%%%%%%%%%%%%%%%%%%%%%%%%%%%%%%%%%%%%%%%%%%%%%%5
\subsection*{The Krichever morphism}

Let $\pi:Y\to X$ be a finite morphism between proper curves over
$\C$. Let us suppose $Y$ and $X$ to be reduced. Fix a set of
pairwise distinct smooth points in $X$, $x=\{x_1,\dots, x_r\}$ and
let $y:=\pi^{-1}(x_1)+\dots+ \pi^{-1}(x_r)$.

Define $A:=H^0(X-x,\o_X)$, $B:=H^0(Y-y,\o_Y)$, $\Sigma_X$ (resp.
$\Sigma_Y$) to be the total quotient ring of $A$ (resp. $B$). Let
$\tr_{\Sigma_X}^{\Sigma_Y}$ denote the trace of $\Sigma_Y$ as a
finite $\Sigma_X$-algebra.

{\it The triple $(Y,X,x)$ is said to have the property $(*)$ if
$\tr_{\Sigma_X}^{\Sigma_Y}(B)\subseteq A$.}

It worth pointing out that every covering $\pi: Y\to X$ has the
property $(*)$ whenever $X$ is smooth or $\pi$ is flat.

Let us fix a set of numerical data as in the previous section
    $$
    E\,=\, \{\bar e_1,\dots, \bar e_r\}
    $$
with $\bar e_i=\{e_1^{(i)},\dots, e_{k_i}^{(i)}\}$ with
$n=e_1^{(i)}+\dots+ e_{k_i}^{(i)}>0$. Let $V$ and $W$ be the
$\C(\!(z)\!)$-algebras defined by the data $E$ as
in~\S\ref{sect:1}. For a $\C$-scheme $S$, we write $\w
V_S:=V\w\otimes_{\C}\o_S$ and $\w W_S:=W\w\otimes_{\C}\o_S$.

\begin{defn}\label{defn:hurwitz}
The Hurwitz functor $\underline{\overline\H}_E^\infty$ of pointed
coverings of curves of degree $n$ with fibres of type $E$ and
formal parameters along the fibers is the contravariant functor on
the category of $\C$-schemes
    {\small $$
    \begin{aligned}
    \underline{\overline\H}_E^\infty \colon
    \left\{{\scriptsize\begin{gathered}\text{category of}\\
    \text{$\C$-schemes}\end{gathered}}\right\}
    & \,\longrightarrow\,
    \left\{{\scriptsize\begin{gathered}\text{category of}\\
    \text{sets}\end{gathered}}\right\}
    \\
    S \,&\, \mapsto   \,
    \underline{\overline\H}_E^\infty(S):=\{(Y,X,\pi,\bar x,\bar y,t_{\bar x},t_{\bar y})\}
    \end{aligned}
    $$}
where
\begin{enumerate}
    \item
    $p_Y:Y\to S$ and $p_X:X\to S$ are proper and flat morphisms whose fibres are
    geometrically reduced curves.
    \item
    $\pi:Y\to X$ is a finite morphism of $S$-schemes of degree $n$
    such that its fibres over closed points $s\in S$ have the property
    $(\ast)$.
    \item
    $\bar x=\{x_1,\dots, x_r\}$ is a set of disjoint smooth sections
    of $p_X$ such that the Cartier divisors $x_i(s)$ for $i=1,\dots,r$
    are smooth points of $X_s:=p_X^{-1}(s)$ for all closed points $s\in
    S$.
    \item
    $\bar y=\{\bar y_1,\dots,\bar y_r\}$ and, for each $i$,
    $\bar y_i=\{y^{(i)}_1,\dots, y^{(i)}_{k_i}\}$ is a set of disjoint smooth sections
    of $p_Y$ such that the Cartier divisor $\pi^{-1}(x_i(S))$ is
    $e^{(i)}_1 y^{(i)}_1(S)+\dots+ e^{(i)}_{k_i} y^{(i)}_{k_i}(S)$.
    \item
    For all closed point $s\in S$ and each irreducible component of
    the fibre $X_s$, there is at least one point $x_i(s)$ lying on
    that component.
    \item
    For all closed point $s\in S$ and each irreducible component of
    the fibre $Y_s$, there is at least one point $y^{(i)}_j(s)$ lying on
    that component.
    \item
    $t_{\bar x}$ is a formal parameter along ${\bar x}(S)$,
    $t_{\bar x}\colon \w \o_{X,{\bar x(S)}}\iso \w V_S$, such
    that it induces
    $$t_{x_i}:=(t_{\bar x})_{x_i}\,\colon\,
    \w \o_{X,{x_i(S)}}\,\iso\, \o_S[\![z_i]\!]$$
    for all $i$.
    \item
    $t_{\bar y}=\{t_{\bar y_1},\dots, t_{\bar y_r}\}$ are formal parameters
    along $\bar y_1(S),\dots, \bar y_r(S)$ such that
    $$\pi^*(t_{x_i})_{y^{(i)}_j(S)}\,=\, t_{y^{(i)}_j}^{e^{(i)}_j}$$
    \item
    $(Y,X,\pi,\bar x,\bar y,t_{\bar x},t_{\bar y})$ and
    $(Y',X',\pi',\bar x',\bar y',t_{\bar x'},t_{\bar y'})$ are said to be
    equivalent when there is a commutative diagram of $S$-schemes
        {\small
        $$\xymatrix{ Y \ar[d]_{\pi} \ar[r]^{\sim} & Y' \ar[d]^{\pi'} \\
        X \ar[r]^{\sim} & X'}$$}
    compatible with all the data.
\end{enumerate}
\end{defn}

The Krichever morphism for the Hurwitz functor is the morphism of
functors
    \beq\label{eq:kricheverforhurwitz}
    \kr\colon \underline{\overline\H}_E^\infty \,\longrightarrow\,
    \underline\gr(W)
    \end{equation}
which sends $(Y,X,\pi,\bar x,\bar y,t_{\bar x},t_{\bar y}) \in
\underline{\overline\H}_E^\infty(S)$ to the following submodule of
$W\w\otimes_{\C}\o_S$
    $$
    t_{\bar y}\Big(\limil{m}
    (p_{Y})_*\o_Y(m\cdot \pi^{-1}(\bar x))\Big)\,\subset\, W\w\otimes_{\C} \o_S
    $$

Let $\M^\infty(\bar r)$ be the moduli functor parametrizing the
classes of sets of data $(Y; y_1,\dots, y_{\bar r}; t_1,\dots,
t_{\bar r})$ of geometrically reduced curves with $\bar r$
pairwise distinct marked smooth  points $\{y_1,\dots, y_{\bar
r}\}$ and formal parameters $\{t_1,\dots, t_{\bar r}\}$ at these
points such that each irreducible component of $Y$ contains at
least one of the marked points. Applying Theorem~4.3 of
\cite{Hurwitz}, one gets the following

\begin{thm}
The Krichever morphism~(\ref{eq:kricheverforhurwitz}) identifies
$\M^\infty(\bar r)(S)$ with the set of submodules $U\in\gr(W)(S)$
such that $U\cdot U\subseteq U$ and $\o_S\subseteq U$. In
particular, this functor is representable by a closed subscheme of
the infinite Grassmannian which will be also denoted by
$\M^\infty(\bar r)$.
\end{thm}

Let $\tr:W\to V$ denotes the trace map of $W$ as a $V$-algebra,
then it is clear that $\tr=\oplus_{i=1}^r\tr^i$ where
$\tr^i:W^i\to \C(\!(z_i)\!)$ is the trace map of the
$\C(\!(z_i)\!)$-algebra $W^i$. Furthermore, one has a commutative
diagram
    $$
    \xymatrix{
    H^0(Y-y,\o_Y) \ar@{^(->}[r]^(.33){t_y}
    \ar[d]_{\tr_{\Sigma_X}^{\Sigma_Y}} &
    \w W_K := W\w\otimes_{\C} K = (W^1\times\dots\times W^r)\w\otimes_{\C} K
    \ar[d]^{\tr}
    \\
    H^0(X-x,\o_X)  \ar@{^(->}[r]^(.33){t_x} &
    \w V_K := V\w\otimes_{\C} K = K(\!(z_1)\!)\times\dots\times K(\!(z_r)\!)}
    $$
for every geometric point $(Y,X,\pi,\bar x,\bar y,t_{\bar
x},t_{\bar y})$ in $\underline{\overline\H}^\infty_E(K)$ ($K$
being an extension of $\C$).

The study of the functor $\underline{\overline\H}^\infty_E$ for
the case $r=1$ has been carried out exhaustively
in~\cite{Hurwitz}. Its main results (given in its fourth section)
are easily generalized for our case. Let us simply state these
generalizations.

\begin{prop}
Let ${\mathcal Y}= (Y,X,\pi,x,y,t_x,t_y)$ be an $S$-valued point
of $\underline{\overline\H}^\infty_E$. Then, it holds that
    $$
    \kr(X,x,t_x)\,=\, \kr({\mathcal Y})\cap \w V_S \,=\, \tr(\kr({\mathcal
    Y}))\,\in \, \gr(V)(S)
    $$
In particular, the Krichever map~(\ref{eq:kricheverforhurwitz}) is
injective.
\end{prop}

\begin{thm}\label{thm:characterizationHwithTr}
Let $U\in \M^\infty(\bar r)\subset \gr(W)$ be an $S$-valued point.
Then, the following conditions are equivalent
\begin{enumerate}
    \item $U\in \underline{\overline\H}^\infty_E(S)$,
    \item $\tr(U)\subseteq U$.
\end{enumerate}
In particular, the functor $\underline{\overline\H}^\infty_E$ is
representable by a closed subscheme ${\overline\H}^\infty_E$ of
$\gr(W)$.
\end{thm}

\begin{thm}[$E$-KP hierarchy]\label{thm:hurwitzbilinear}
Let $B\in\M^{\infty}(\bar r)\subset \gr^m(W)$ be a closed point.
Let $\{u_{1,1},\dots, u_{r,k_r}\}$ be integer numbers defined by
$v_m=z_1^{u_{1,1}/e_1^{(1)}}\dots z_r^{u_{r,k_r}/e_{k_r}^{(r)}}$.

Then, $B\in \overline\H^{\infty}_E$ if and only if the following
``bilinear identities'' are satisfied; that is, the form
    {\small
    $$
    \Big(\sum_{l=1}^r\sum_{j=1}^{k_l} \sum_{i=1}^{e_j^{(l)}}
    \frac{\psi_{a,b,B}^{(l,j)}(\xi_{e_j^{(l)}}^i z_l^{1/{e_j^{(l)}}},t)}
    {(\xi_{e_j^{(l)}}^i z_l^{1/{e_j^{(l)}}})^{\delta_{al}\delta_{bj}-u_{l,j}}}
    \Big)
    \Big(
    \sum_{l=1}^r\sum_{j=1}^{k_l} \sum_{i=1}^{e_j^{(l)}}
    \frac{\psi_{c,d,B}^{(l,j)}(\xi_{e_j^{(l)}}^i z_l^{1/{e_j^{(l)}}},t)}
    {(\xi_{e_j^{(l)}}^i z_l^{1/{e_j^{(l)}}})^{u_{l,j}-1+\delta_{cl}\delta_{dj}}}
    \Big)\frac{dz}z
    $$}
has residue zero at $z=0$ for all $a,b,c,d$ and where $\xi_e$ is a
primitive $e$-th root of unity and $\delta$ is the Kronecker
symbol.
\end{thm}

\begin{pf}
The proof is modeled in the proof of Theorem~4.10 or
\cite{Hurwitz}. The idea consists of translating the second
condition of the previous statement in terms of Baker-Akhiezer
functions by means of Theorem~\ref{thm:BAgeneragrmV}.
\end{pf}

\begin{defn}
The Hurwitz functor $\underline{\H}_E^\infty$ of pointed coverings
of smooth curves of degree $n$ with fibres of type $E$ and formal
parameters along the fibers is the subfunctor of
$\underline{\overline\H}_E^\infty$ consisting of data
$(Y,X,\pi,\bar x,\bar y,t_{\bar x},t_{\bar y})$ where the fibres
of $Y\to S$ are nonsingular curves.
\end{defn}

\begin{thm}
The functor $\underline{\H}_E^\infty$ is representable by a
subscheme of $\gr(W)$ which will be denoted by ${\H}_E^\infty$.
\end{thm}

\begin{pf}
Observe that if $Y\to X$ is a family has the property $(*)$ and
the fibres of $Y$ are nonsingular curves, then the fibres of $X$
are also nonsingular curves. Now let ${\mathcal Y}\to
{\overline\H}^\infty_E$ be the universal family given by the
representability of the functor
$\underline{\overline\H}^\infty_E$. Then the desired subscheme,
${\H}_E^\infty$, consists precisely of the points $s\in
{\overline\H}^\infty_E$ such that ${\mathcal Y}_s$ is smooth.
\end{pf}

\begin{rem}
Standard procedures can be used to rewrite the identities of
Theorems~\ref{thm:bilinearGrW} and~\ref{thm:hurwitzbilinear} as
hierarchies of partial differential equations for the
$\tau$-function.
\end{rem}

%%%%%%%%%%%%%%%%%%%%%%%%%%%%%%%%%%%%%%%%%%%%%%%%%%%%%%%%%%%%%%%%%%%%%5
\subsection*{Coverings with prescribed ramification}

The above results, which generalize the results proved
in~\cite{Hurwitz}, allow us to describe the structure of the
Hurwitz scheme parametrizing coverings with prescribed
ramification points and formal parameters at the given points.

Let $\pi: Y\to X$ be a finite covering of degree $n$ of smooth
integral curves. Let $\bar g$ and $g$ be the genus of $Y$ and $X$
respectively. Then, the Hurwitz formula reads
    $$
    1-g\,=\, n(1-g)-\frac12 \sum_{y\in Y} (e_y-1)
    $$
where $e_y$ is the ramification index of the point $y\in Y$.

Let $\{x_1,\dots,x_r\}\subset X$ be the branch locus and denote
    $$
    \begin{gathered}
    \bar y\,=\, \pi^{-1}(x_1)+ \dots + \pi^{-1}(x_r)
    \\
    \pi^{-1}(x_i)\,=\, e_1^{(i)}y_1^{(i)}+\dots + e_{k_i}^{(i)}y_{k_i}^{(i)}
    \end{gathered}
    $$
Considering $E= \{\bar e_1,\dots, \bar e_r\}$, $\bar e_i=
\{e_1^{(i)}, \dots , e_{k_i}^{(i)}\}$ and $\bar r= \sum_{i=1}^r
k_i$, one has that
    $$
    \sum_{y\in Y} (e_y-1)\,=\,  \sum_{i=1}^r \sum_{j=1}^{k_i}(e_j^{(i)}-1)
    \,=\, r n -\bar r
    $$
then the Hurwitz formula can be rewritten as
    $$
    1-g \,=\, n(1-g)-\frac12(rn-\bar r)
    $$

\begin{defn}
For integers $i,j$, we define the following subschemes of
${\H}^\infty_E$
    $$
    \begin{gathered}
    {\H}^\infty_E[j]\,:=\,\{
    U\in {\H}^\infty_E\cap \gr^{1-j}(W)\text{ such that
    $U$ is an integral domain}\}
    \\
    {\H}^\infty_E[j,i]\, :=\,
    \{U\in {\H}^\infty_E[j]\text{ such that }\tr(U)\in
    \gr^{1-i}(V)\}
    \end{gathered}
    $$
\end{defn}

Note that, since $\tr(U)\subseteq U$ for any $U\in {\H}^\infty_E$,
the condition that $U$ is an integral domain implies that $\tr(U)$
is integral too.

From the representability of ${\H}^\infty_E$ and the Hurwitz
formula for coverings, one has the following

\begin{thm}\label{thm:Hurwitz[j,i]}
Let $(E,n,\bar r)$ be a set of numerical data as in~\S\ref{sect:1}
and $g,\bar g$ be two non-negative integer numbers  satisfying
    $$
    \bar g -1\,=\, n(g-1) +\frac12(rn-\bar r)
    $$

Then, the subscheme ${\H}^\infty_E[\bar g,g]\subset \gr^{1-\bar
g}(W)$ is the moduli scheme parametrizing geometrical data
$(Y,X,\pi,\bar x,\bar y,t_{\bar x},t_{\bar y})\in {\H}^\infty_E$
where $Y$ has genus $\bar g$, $X$ has genus $g$, the covering
$\pi:Y\to X$ is non-ramified outside $\bar y$ and the ramification
index at the point $y_j^{i}\in Y$ is $e_j^{(i)}$.
\end{thm}

\begin{rem}
Let us observe that there exists a natural forgetful morphism
    $$
    \begin{aligned}
    \Phi\colon {\H}^\infty_E[\bar g,g] \,& \longrightarrow\,
    \M^\infty_g(r)
    \\
    (Y,X,\pi,\bar x,\bar y,t_{\bar x},t_{\bar y}) &\longmapsto (X,{\bar x},t_{\bar x})
    \end{aligned}
    $$
Given $(X,{\bar x},t_{\bar x})\in \M^\infty_g(r)$ with $X$ smooth
and integral, the fiber of this point will be denoted by
${\H}^\infty_E(X,\bar x,t_{\bar x})$. Recall that there is a
finite number of coverings $\pi:Y\to X$ with $Y$ and $X$ smooth,
$X$ integral, $\pi$ is dominant on each component of $Y$,
non-ramified outside $\bar x$, and ramification indexes $\bar e_i$
at $x_i$ (see \cite{elsv,OP}). Then, we conclude that
${\H}^\infty_E(X,\bar x,t_{\bar x})$ is a finite set.
\end{rem}

%%%%%%%%%%%%%%%%%%%%%%%%%%%%%%%%%%%%%%%%
\section{Tangent space to the Hurwitz scheme}
%%%%%%%%%%%%%%%%%%%%%%%%%%%%%%%%%%%%%%%%

This section is devoted to an explicit computation of the tangent
space of the Hurwitz schemes constructed in the previous section.
To this goal, we begin by recalling the computation of the tangent
spaces to the infinite Grassmannian and to the moduli space of
pointed curves.

\begin{prop}\label{prop:tangentgr}
Let $U$ be a rational point of $\gr(W)$.  There is a canonical
isomorphism
    $$
    T_U\gr(W)\,\overset{\sim}\longrightarrow \, \hom(U,W/U)
    $$
\end{prop}

\begin{pf}
Let $A\sim W_+$ and $U\in\gr(W)$ such that $U\oplus A\simeq W$.
Let $F_A$ be the open subscheme of $\gr(W)$ parametrizing those
subspaces $U'\in\gr(W)$ such that $U'\oplus A\simeq W$. Then, it
is well known that $F_A$ is isomorphic to the affine space
$\hom(U,A)$ and that the embedding
    $$
    \hom(U,A)\,\iso\, F_A
    \, \hookrightarrow\, \gr(W)
    $$
maps $f:U\to A$ to its graph $\Gamma_f:=\{u+f(u) \vert u\in U\}$.

Since $U\in F_A$ (it corresponds to the zero map) and $F_A$ is
open, we obtain an isomorphism of vector spaces
    {\small $$
    \hom(U,A) \,\overset{\sim}\to\, T_0 \hom(U,A)
    % \,\simeq \,
    % F_A(k[\epsilon]/\epsilon^2)\times_{F_A(k)}\{0\}
    \,\overset{\sim}\to\,
    T_U\gr(W) \,\overset{\sim}\to\,
    \gr(W)(k[\epsilon]/\epsilon^2)\times_{\gr(W)(k)}\{ U \}
    $$}
which maps $f\in \hom(U,A)$ to the
$(k[\epsilon]/\epsilon^2)$-valued point of $\gr(W)$ given by
$\{u+\epsilon f(u) \vert u\in U\}$.

Composing the inverse of this map with the isomorphism
    $$
    \begin{aligned}
    \hom(U,A) \,& \overset{\sim}\longrightarrow\, \hom(U,W/U)
    \\
    f\,& \longmapsto \quad{\small  U\overset{\id+f}\to U\oplus A\overset{\sim}\to W\to
    W/U}
    \end{aligned}
    $$
one gets the desired isomorphism (observe that it does not depend
on the choice of $A$).
\end{pf}

\begin{prop}\label{prop:tangentM}
Let $U$ be a rational point of $\M^\infty(\bar r)$.  The
isomorphism of Proposition~\ref{prop:tangentgr} induces a
canonical identification
    $$
    T_U\M^\infty(\bar r)\,\simeq \, \der(U,W/U)
    $$
(where $\der$ means derivations trivial over $\C$).

Furthermore, if $U$ is associated to the geometrical data $(C,\bar
p,\bar z)$ under the Krichever map, then $W/U\simeq H^0(C-\bar
p,\omega_C)^*$.
\end{prop}

\begin{pf}
For $U\in\M^\infty(\bar r)$, one has that
    $$
    T_U \M^\infty(\bar r)\,=\,
    \{\bar U\in T_U\gr(W)\text{ such that } \bar U\cdot\bar U\subseteq
    \bar U\text{ y } k[\epsilon]\subset \bar U\}
    $$
From Proposition~\ref{prop:tangentgr}, there is a map $f\in
\hom(U,A)$ such that $\bar U=\{u+\epsilon f(u) \vert u\in U\}$.

The condition $\bar U\cdot\bar U\subseteq \bar U$ means that for
$u,u'\in U$ there exists $u''\in U$ satisfying $(u+\epsilon
f(u))\cdot (u'+\epsilon f(u'))= u''+\epsilon f(u'')$; that is
    \beq\label{eq:derivacion1}
    f(u\cdot u')\,=\, u f(u')+f(u) u'
    \end{equation}
The second condition, $k[\epsilon]\subset \bar U$, implies that
there exists $u_0\in U$ such that $u_0+\epsilon f(u_0)=1$; or, in
other words
    \beq\label{eq:derivacion2}
    f(1)\,=\,0
    \end{equation}
It is now easy to check that the image of $f\in\hom(U,A)$ in
$\hom(U,W/U)$ gives rise to a derivation $D_f\in\der(U,W/U)$ (note
that $W/U$ is an $U$-module). The claim follows from a
straightforward check.

The second part of the statement follows easily from the exact
sequence
    {\small $$
    \xymatrix@=12pt{
    0\ar[r] &  H^0(C-\bar p,\o_C)\ar[r]\ar[d]^{\wr}  &
    W\simeq (\w \o_{\bar p})_{(0)} \ar[r]\ar[d]^{\wr}  &
    \limpl{m} H^1(C,\o_C(-m\bar p)) \ar[r] & 0
    \\
    & U & W
    }$$}
%
%    {\small $$
%    0\, \to\, U\simeq H^0(C-\bar p,\o_C) \,\to\,
%    W\simeq (\w \o_{\bar p})_{(0)}\,\to \, \limpl{m}
%    H^1(C,\o_C(-m\bar p))\,\to\,0
%    $$}
\end{pf}

Note that a similar result holds for points of
$\M^\infty(r)\subset \grv$.

Let us denote
    {\small
    \beq\label{eq:derTr}
    \der(U,W/U)^{\tr}\,:=\, \left\{
    \text{$D\in\der(U,W/U)$  such that
    \raisebox{25pt}{
    \xymatrix{ U \ar[r]^D \ar[d]_{\tr} & W/U \ar[d]^{\tr^1} \\
    \tr U \ar[r]^D & V/\tr U }}}
    \right\}
    \end{equation}}
where $\tr^1$ is the map induced by the trace (since $\tr
U\subseteq U$).

\begin{thm}\label{thm:tangentehurwitz}
Let $U$ be a rational point of ${\overline\H}^\infty_E$. The
embedding ${\overline\H}^\infty_E\hookrightarrow\M^\infty(\bar r)$
yields an identification
    $$
    T_U{\overline\H}^\infty_E\,\simeq \,
    \der(U,W/U)^{\tr}
    $$

Moreover, if $U$ corresponds to the geometrical data
$(Y,X,\pi,\bar x,\bar y,t_{\bar x},t_{\bar y})$ in $\H^\infty_E$,
then $\tr^1$ is the map
    $$
    \tr^1\colon H^0(Y-y,\omega_Y)^*\,\longrightarrow\,
    H^0(X-x,\omega_X)^*
    $$
canonically induced by the trace $\pi_*\o_Y\to\o_X$.
\end{thm}

\begin{pf}
Let us keep the notations of the proofs of the two previous
propositions.  Then, let $\bar U\in\gr(W)$ be a
$k[\epsilon]/\epsilon^2$-valued point lying in $T_U\M^\infty(\bar
r)$. Let $f\in\hom(U,A)$ correspond to $\bar U$.

Then, the condition $\tr(\bar U)\subseteq \bar U$
(Theorem~\ref{thm:characterizationHwithTr}) is equivalent to say
that for each $u\in U$ there exists $u'\in U$ satisfying
$\tr(u+\epsilon f(u))= u'+\epsilon f(u')$. Since $\tr(u)\in U$,
this condition is
    \beq\label{eq:ftr}
    \tr^1(f(u))\,=\, f(\tr(u)) \qquad \forall u\in U
    \end{equation}
And the first part of the claim follows easily.

Let us prove the second part. Note that the trace is a sheaf
homomorphism
    $$
    \tr\colon \pi_*\o_Y\,\longrightarrow\,\o_X
    $$
This map induces
    $$
    H^1(X, (\pi_*\o_Y)(-n\bar x))\, \longrightarrow
    \, H^1(X, \o_X(-n \bar x))
    $$
Bearing in mind the adjunction formula, that $\pi$ is affine,
Serre duality and taking limits, one obtains the desired map
    $$
    \tr^1\colon H^0(Y-\bar y,\omega_Y)^*\,\longrightarrow\,
    H^0(X-\bar x,\omega_X)^*
    $$
Using arguments similar to those of the proof of
Proposition~\ref{prop:tangentM}, it is not difficult to check that
this map is compatible with the isomorphisms $H^0(Y-\bar
y,\omega_Y)^*\simeq W/U$ and $H^0(X-\bar x,\omega_X)^*\simeq
V/\tr(U)$.
\end{pf}

\begin{thm}\label{thm:tangenteHurwitzinjectivo}
Let $(E,n,r,\bar g,g)$ as in Theorem~\ref{thm:Hurwitz[j,i]}. Let
$U\in \H^\infty_E[\bar g,g]$ be a rational point.
% corresponding to data $(Y,X,\pi,\bar x,\bar y,t_{\bar x}, t_{\bar y})$.
Then, there is a canonical injection
    $$
    T_U \H^\infty_E[\bar g,g] \,\hookrightarrow\, T_{\tr U} \M^\infty(r)
    $$
\end{thm}

\begin{pf}
By Proposition~\ref{prop:tangentM} and
Theorem~\ref{thm:tangentehurwitz}, the claim is equivalent to the
injectivity of the restriction map
    $$
    \der(U,W/U)^{\tr} \,\longrightarrow\, \der(\tr U,V/\tr U)
    $$

Let  $D\in\der(U,W/U)^{\tr}$ be in the kernel of the above map;
that is, $D\vert_{\tr(U)}=0$. Let $u$ be any element in $U$ and
let us see that $Du=0$. Since $A:=\tr(U)\to U$ is an integral
morphism, there is a monic minimal polynomial $p(x)=\sum_i a_i
x^i\in A[x]$ such that $p(u)=0$. Then, the following identity
holds true
    $$
    0\,=\, D p(u) \,=\,
    \sum_i (D a_i)u^i + p'(u) Du \,=\,
    p'(u) Du
    $$
($p'(x)$ denoting the derivative w.r.t. $x$).

Since $p(x)$ is separable and $\pi$ is unramified over $X-\bar
x=\sp A$, then $\frac{1}{p'(u)}\in A[u]\subseteq U$. Therefore,
one has $Du=0$.
\end{pf}

%%%%%%%%%%%%%%%%%%%%%%%%%%%%%%%%%%%%%%%%
\section{The multicomponent Virasoro group}
%%%%%%%%%%%%%%%%%%%%%%%%%%%%%%%%%%%%%%%%

The algebraic Virasoro group, defined as
$G:=\underline\aut_{\C}\C(\!(z)\!)$, was introduced and studied
in~\cite{automorphism}. In this section it will be generalized for
certain $\C(\!(z)\!)$-algebras. Note that $V$ carries the linear
topology given by $\{z^n V_+\}_{n\in \Z}$.

It is convenient to recall the following result: {\sl let $R$ be a
$\C$-algebra and let $f(z)\in R(\!(z)\!)$. Then $f(z)$ is
invertible if and only if there exists $n\in\Z$ such that
$a_i\in\rad(R)$ for $i< n$ and $a_n$ is invertible}.

\begin{defn}
The functor of automorphisms of $V$ is the functor defined from
the category of $\C$-schemes to the category of groups defined as
follows
    $$
    S\,\rightsquigarrow\,
    \underline\aut_{\C}(V)(S):= \aut_{H^0(S,\o_S)}
    V\w\otimes_{\C} H^0(S,\o_S)
    $$
where $\aut_{R}$ means continuous automorphisms of $R$-algebras.
\end{defn}

\begin{lem}
Let $S$ be a $\C$-scheme and let $\phi\in
\underline\aut_{\C}(V)(S)$. For any point $s\in S$, let
$p_{\phi(s)}$ be the permutation defined by $\phi(s)$ on the set
$\sp(V_s)$ (note that $V_s:=V\w\otimes_{\C} k(s)$ consists of $r$
points).

Then, the map
    $$
    \begin{aligned}
    S\,&\longrightarrow\, {\mathcal S}_r \\
    s\,&\longmapsto \, p_{\phi(s)}
    \end{aligned}
    $$
(${\mathcal S}_r$ being the symmetric group of $r$ letters) is
locally constant.
\end{lem}

\begin{pf}
Let us assume that $S=\sp(R)$ is an irreducible $\C$-scheme. Let
$\phi\in \underline\aut_{\C}(V)(S)$. Since $V_s:=V\w\otimes_{\C}
k(s)\simeq \prod^r k(s)(\!(z)\!)$ is a product of fields, it
follows that $\phi(s)\in \underline\aut_{\C}(V)(k(s))$ acts by
permutation on $\sp(V_R)$; that is, on the set of ideals
$I_1:=((z_1,0,\dots, 0))$, \dots, $I_r:=((0,\dots, 0,z_r))$. That
is, there exists $p_{\phi(s)}\in{\mathcal S}_r$ such that
    $$
    \phi(s)(I_i)\,=\, I_{p_{\phi(s)}(i)}
    $$

Let us write $\phi(0,\dots,z_i,\dots,0)=(\phi_1^{(i)},\dots
,\phi_r^{(i)})\in R((z_1))\times\dots\times R((z_r))$. Let $s_0$
denote the point associated to the minimal prime ideal. The very
definition of the permutation $p_{\phi(s_0)}$ says that
    $$
    \phi_j^{(i)}(s_0)\in k(s_0)((z_j)) \text{ is }
    \begin{cases}
    & \text{invertible if }j=p_{\phi(s_0)}(i)
    \\
    & \text{zero if }j\neq p_{\phi(s_0)}(i)
    \end{cases}
    $$
and, therefore, $s_0\in (\phi_j^{(i)})_0$ for $j\neq
p_{\phi(s_0)}(i)$. Since $s_0$ is minimal, $S$ is irreducible and
$(\phi_j^{(i)})_0$ is closed, it follows that the closure of
$s_0$, $S$, is contained in $(\phi_j^{(i)})_0$  for $j\neq
p_{\phi(s_0)}(i)$.

Then, it turns out that $\phi_j^{(i)}(s)\neq 0$  for $j=
p_{\phi(s_0)}(i)$ and for all $s\in S$; equivalently,
$\phi_j^{(i)}\in R((z_j))$ is invertible for $j= p_{\phi(s_0)}(i)$
and, therefore, $p_{\phi(s_0)}=p_{\phi(s)}$ for all $s\in S$. The
statement follows.
\end{pf}

As a consequence of the previous lemma we have a well defined map
of functors on groups
    $$
    \underline\aut_{\C}(V)\,\overset{p}\longrightarrow\, {\mathcal S}_r
    $$
given by $p(\phi):=p_\phi$. Here, the scheme structure of the
finite group ${\mathcal S}_r$ is considered to be isomorphic to
the $\C$-scheme $\sp\big(\overset{r!}\prod\C\big)$.

\begin{thm}
Let $V$ be the  $\C$-algebra $\C(\!(z_1)\!)\times \dots\times
\C(\!(z_r)\!)$ and $V_+$ the subalgebra $\C[\![z_1]\!]\times
\dots\times \C[\![z_r]\!]$.

The canonical exact sequence of functors on groups
    $$
    0\to G\times \overset{r}\dots\times G
    \overset{i}\longrightarrow
    \underline\aut_{\C}(V)\overset{p}\longrightarrow {\mathcal S}_r
    \to 0
    $$
splits.

In particular, $\underline\aut_{\C}(V)$ is representable by a
formal group $\C$-scheme which will be denoted by $G_V$.
\end{thm}

\begin{pf}
The map $i$ is the canonical inclusion
    $$
    \underline\aut_{\C}\C(\!(z_1)\!)\times\dots\times \underline\aut_{\C}\C(\!(z_r)\!)
    \,\hookrightarrow\,
    \underline\aut_{\C}(V)
    $$

Now, observe that one can associate to every permutation
$\sigma\in {\mathcal S}_r$ the automorphism of $V$ defined by
$z_i\mapsto z_{\sigma^{-1}(i)}$.
\end{pf}

\begin{cor}
If $G_V^0$ denotes the connected component of the origin of $G_V$,
then there is an identification of groups
    $$
    G_V^0\,=\, G \times\overset{r} \dots\times G
    $$

There are scheme isomorphisms
    $$
    \begin{gathered}
    G_V\,\simeq \, (G\times \overset{r}\dots\times G)\times {\mathcal
    S}_r
    \\
    G_V^0\,\simeq \, \Gamma_V
    \end{gathered}
    $$
(which are not group homomorphisms).
\end{cor}

The formal group scheme $G_V^0$ acts on $\gr(V)$ and this action
yields an action on the subscheme $\M^\infty(r)\subset\gr(V)$.

The main results of~\cite{automorphism} can be generalized to the
present setting.

Now, let us consider the groups $G_V$ and $G_W$ corresponding to
the $\C$-algebras $V$ and $W$ respectively.

The $V$-algebra structure of $W$ allows us to define a subgroup
$G_V^W\subseteq G_W^0$ as follows
    \beq\label{defn:gvw}
    G_V^W\,:=\,
    \left\{
    \raisebox{20pt}{\xymatrix{
        W\ar[r]_{\sim}^{\bar g} & W \\
        V \ar@{^(->}[u] \ar[r]_{\sim}^{g} & V \ar@{^(->}[u]
        }}
    \quad\text{where $\bar g\in G_W^0$ and $g\in G_V^0$}
    \right\}
    \end{equation}
The element of $G_V^W$ given by such a diagram will be denoted by
$\bar g$.

\begin{prop}
There exists an exact sequence of formal group $\C$-schemes
    $$
    0 \to \mu_E \to G_V^W \overset{\pi}\to G_V^0 \to 0
    $$
where $\pi(\bar g)=g$ and
    $$
    \mu_E \,:=\,
    (\mu_{e_1^{(1)}}\times \dots\times\mu_{e_{k_1}^{(1)}})
    \times \dots\times
    (\mu_{e_1^{(r)}}\times \dots\times\mu_{e_{k_r}^{(r)}})
    $$
($\mu_{e_i^{(j)}}\subset \C^*$ being the group of $e_i^{(j)}$-th
roots of unity).
\end{prop}

\begin{pf}
By the previous Corollary it can be assumed that $r=1$ and $k_1=$;
that is, $V=\C(\!(z)\!)$ and $W=\C(\!(z^{1/e})\!)$. Then, one has
that $G_V\simeq G$ and $G_W\simeq G$.

Observe that an element $\bar g\in G_{W}^0$ belongs to $G_V^W$ if
and only if $g(z)= \bar g(z^{1/e})^{e}$. And the result follows.
\end{pf}

\begin{cor}\label{cor:LieGVWisoLieGV}
The canonical restriction map $G_V^W\to G_V^0$ yields an
isomorphism of Lie algebras
    $$
    \lie G_V^W\,\overset{\sim}\longrightarrow\, \lie G_V^0
    $$
\end{cor}

\begin{lem}
Let $R$ be a $\C$-algebra and $f(z^{1/e}) \in R(\!(z^{1/e})\!)$.

If $f(z^{1/e})^e \in R(\!(z)\!)$ and $f(z^{1/e})$ is invertible,
then there exist $i$ such that $z^{i/e}f(z^{1/e}) \in R(\!(z)\!)$.
\end{lem}

\begin{pf}
We may assume that $f(z^{1/e})=\sum_{i}a_i z^{i/e}$ where $a_i\in
\rad(R)$ for $i< 0$ and $a_0$ is invertible. Let $I\subseteq
\rad(R)$ be the ideal generated by $\{a_i \vert i< 0\}$ (recall
that this set is finite) and let $n\geq 0$ be the smallest integer
such that $I^{n+1}$ vanishes. Let us proceed by induction on $n$.

Case $n=1$; that is $a_i\cdot a_j=0$ for all negative integers
$i,j$. Let $i_0$ be the smallest index for which $a_{i_0}\neq 0$.
The hypothesis implies that
    $$
    f(z^{1/e})^e \,=\,
    e\cdot a_0^{e-1}\cdot a_{i_0}\cdot z^{i_0/e}+ (\text{higher order terms})
    \,\in\,  R(\!(z)\!)
    $$
and, therefore, $i_0=\dot{e}$. To conclude is suffices to consider
 $(a_0+ a_{i_0}z^{i_0/e})^{-1}\cdot f(z^{1/e})$ and iterate
this argument.

General case. Let $\bar f_n(z^{1/e})$ be the class of $f(z^{1/e})$
in $R/I^n(\!(z^{1/e})\!)$. From the induction hypothesis it
follows that $\bar f_n(z^{1/e})\in R/I^n(\!(z)\!)$. Let $f_n(z)\in
R(\!(z)\!)$ be a preimage of $\bar f_n$. Consider now the element
element $\frac{f(z^{1/e})}{f_n(z)} \in R(\!(z^{1/e})\!)$. From the
$n=1$ case, it follows that $f(z^{1/e}) (f_n(z))^{-1} \in
R(\!(z)\!)$, and the claim follows.
\end{pf}

\begin{lem}
Let $R$ be a $\C$ algebra, $\tr: \w W_R\to \w V_R$ be $\w V_R$ be
the trace map and $\bar g$ be an element of $G^0_W$.

Then, $\bar g\in G_V^W$ if and only if $\tr\circ \bar g= \bar
g\circ \tr$.
\end{lem}

\begin{pf}
Note that it suffices to prove the claim for the following case
    $$
    V=R(\!(z)\!)\,\hookrightarrow\, W=R(\!(z^{1/e})\!)
    $$
Let us show that
    \beq\label{eq:identitytrace}
    \tr({\bar g} w)\,=\, \pi({\bar
    g})\tr(w)
    \,=\, g(\tr(w))\qquad \forall w\in\w W_R
    \end{equation}
An element $\bar g\in G_W^0(R)$ is of the form
    $$
    \bar g(z^{1/e})\,=\, z^{1/e}\cdot (\sum_i a_i z^{i/e})\,\in\,
    R(\!(z^{1/e})\!)
    $$
where $a_i$ is nilpotent for $i<0$ and $a_0$ in invertible. The
condition $\bar g\in G_V^W$ implies that
    $$
    \bar g(z^{1/e})^e \,=\, \bar g(z) \,\in\, R(\!(z)\!)
    $$
and, by the previous lemma, it follows that $\sum_i a_i z^{i/e}\in
R(\!(z)\!)$, that is, $a_i=0$ if $i\neq \dot{e}$. In particular,
$g= \pi(\bar g)$ is the automorphism given by
    $$
    g(z)\,=\, z (\sum_j a_{j e} z^j)\,\in\, R(\!(z)\!)
    $$
By linearity, it is enough to check the claim for $z^{l/e}$ where
$l\in \Z$
    $$
    \tr(\bar g(z^{l/e}))\,=\,
    \tr(z^{l/e}\cdot (\sum_j a_{j e} z^j))^l)\,=\,
    (\sum_j a_{j e} z^j))^l\cdot \tr(z^{l/e})\,=\,
    g(\tr(z^{l/e}))
    $$
since $\tr(z^{l/e})=0$ for $l\neq \dot{e}$ and
$\tr(z^{l/e})=e\cdot z^{l/e}$ for $l= \dot{e}$.

Conversely. For $\bar g\in G^0_W$ commuting with the trace, we
have that
    $$
    \begin{gathered}
    \bar g(\tr z)\,=\, \bar g(e\cdot z)\,=\,
    e\cdot \bar g(z^{1/e})^e
    \\
    \tr(\bar g(z))\,\in\, R(\!(z)\!)
    \end{gathered}
    $$
Therefore, one obtains that $\bar g(z^{1/e})^e\in  R(\!(z)\!)$.
The previous lemma implies that $\bar g(z^{1/e})$ is of the form
$z^{1/e}(\sum_i a_i z^i)$ which belongs to $G_V^W$.
\end{pf}

\begin{thm}\label{thm:LieGVWisoDerW}
It holds that
    \begin{enumerate}
    \item $\lie G_V\simeq \oplus_{i=1}^r \der(\C(\!(z_i)\!),\C(\!(z_i)\!))
        \simeq \oplus_{i=1}^r \C(\!(z_i)\!){\frac{\partial}{\partial {z_i}}}$;
    \item $\lie G_W\simeq \oplus_{i=1}^r \oplus_{j=1}^{k_i}
        \der(\C(\!(z_i^{1/e_j^{(i)}})\!),\C(\!(z_i^{1/e_j^{(i)}})\!))$;
    \item $\lie G_V^W\simeq (\lie G_W)^{\tr}$ (those derivations commuting with the trace; as in~(\ref{eq:derTr})).
    \end{enumerate}
\end{thm}

\begin{pf}
The first two follow from similar arguments as in the proof of
Theorem~3.5 of \cite{automorphism}. The last claim is a
consequence of the previous Lemma.
\end{pf}

\begin{thm}
The group $G_V^W$ acts on ${\overline\H}^\infty_E$.
\end{thm}

\begin{pf}
Recall that ${\overline\H}^\infty_E$ consists of those points $U$
of $\M^\infty(\bar r)$ such that $\tr(U)\subseteq U$ and that
$G_V^W$ acts on $\M^\infty(\bar r)$. Therefore, it suffices to
show that the condition $\tr(U)\subseteq U$ implies that $\tr(\bar
g U)\subseteq \bar g U$ for all $\bar g\in G_V^W$.

Let $U$ be a point of ${\overline\H}^\infty_E$. From the defining
property of $\bar g$ (see~(\ref{defn:gvw})) and from the inclusion
$\tr(U)\subseteq U$, one has
    $$
    \tr(\bar g U)\,=\, \bar g(\tr(U)) \,\subseteq\, \bar g U
    $$
and the statement follows.
\end{pf}

\begin{thm}\label{thm:GonHloctran}
Let $(E,n,r,\bar g,g)$ as in Theorem~\ref{thm:Hurwitz[j,i]}.
 The group $G_V^W$ acts on
${\H}^\infty_E[\bar g, g]$ and this action is locally transitive.
\end{thm}

\begin{pf}
It is straightforward that the action on ${\overline\H}^\infty_E$
gives an action on ${\H}^\infty_E[\bar g,g]$.  The rest of the
proof is based on the ideas of~\cite{automorphism} (Lemma~4.10,
Theorem~4.11 and Lemma~A.2). There it is shown that it is enough
to check that the surjectivity of the map of tangent spaces
    \beq\label{eq:surjectivitytangetgvwH}
    T_{\id} G_V^W\,\longrightarrow\, T_U \H^\infty_E[\bar g,g]
    \end{equation}

Let $(Y,X,\pi,\bar x,\bar y,t_{\bar x},t_{\bar y})$ be the data
attached to $U\in \H^\infty_E[\bar g,g]$. Then,
Theorems~\ref{thm:tangentehurwitz},
\ref{thm:tangenteHurwitzinjectivo} and~\ref{thm:LieGVWisoDerW} and
Corollary~\ref{cor:LieGVWisoLieGV} give the following commutative
diagram
    {\small
    $$
    \xymatrix{
    0 \ar[r] &  H^0(Y-\bar y,\T_Y) \ar[r] & \lie G_W^0 \ar[r] &
    \der(U,W/U) \ar[r] & 0 \\
    & & \lie G_V^W \ar@{^(->}[u] \ar[r]^-{\psi} \ar[d]^{\wr} &     \der(U,W/U)^{\tr}
    \ar@{^(->}[u] \ar@{^(->}[d] \\
    0 \ar[r] &  H^0(X-\bar x,\T_X) \ar[r] & \lie G_V^0 \ar[r] &
    \der(\tr U,V/\tr U) \ar[r] & 0
    }
    $$}
($\T$ denoting the tangent sheaf). From the diagram, we deduce
that $\psi$ is surjective. Since $\psi$ coincides with the
map~(\ref{eq:surjectivitytangetgvwH}), we are done.
\end{pf}

\begin{thm}
Let $(E,n,r,\bar g,g)$ as in Theorem~\ref{thm:Hurwitz[j,i]}. Let
$U\in \H^\infty_E[\bar g,g]$ be a rational point.
% corresponding to data $(Y,X,\pi,\bar x,\bar y,t_{\bar x}, t_{\bar y})$.
Then, there is an isomorphism
    $$
    T_U \H^\infty_E[\bar g,g] \,\overset{\sim}\longrightarrow\, T_{\tr U} \M^\infty(r)
    $$
\end{thm}

\begin{pf}
The injectivity follows from
Theorem~\ref{thm:tangenteHurwitzinjectivo}. The surjectivity is a
consequence of the diagram of the previous proof.
\end{pf}

\begin{rem}
This Theorem is the analog of the fact that the map from the
classical Hurwitz space to the moduli of curves is etale at those
points corresponding to covers where both curves are smooth. Our
approach also allows one to study the non-smooth case, however,
because of our goals we have focus ourselves on the smooth case.
\end{rem}

%%%%%%%%%%%%%%%%%%%%%%%%%%%%%%%%%%%%%%%%%%
\section{Picard schemes}
%%%%%%%%%%%%%%%%%%%%%%%%%%%%%%%%%%%%%%%%%%

\begin{defn}\label{defn:pic}
Let $\underline\pic^\infty_E[\bar g,g]$ be the contravariant
functor from the category of $\C$-schemes to the category of sets
defined by
    $$
    S\,\rightsquigarrow\, \{(Y,X,\pi,\bar x,\bar y,t_{\bar
    x},t_{\bar y},L,\phi_{\bar y})\}
    $$
where
    \begin{enumerate}
    \item $(Y,X,\pi,\bar x,\bar y,t_{\bar x},t_{\bar y})\in \H_E^\infty[\bar
    g,g](S)$;
    \item $L$ is a line bundle on $Y$ such that $(Y_s,\bar y(s),L\vert_{Y_s})$ is
    maximal for all closed point $s\in S$;
    \item $\phi_{\bar y}$ is a formal trivialization of $L$
    along $\bar y$; that is, an isomorphism $\phi_{\bar y} : \w
    L_{\bar y}\simeq \w\o_{Y,\bar y}$.
    \item two sets of data are said to be
    equivalent when there is an isomorphism of $S$-schemes
    $Y\iso Y'$ compatible with all the data.
    \end{enumerate}
\end{defn}

\begin{thm}
The functor $\underline\pic^\infty_E[\bar g,g]$ is representable
by a subscheme $\pic^\infty_E[\bar g,g]$ of $\gr(W)$.
\end{thm}

\begin{pf}
Consider the morphism from $\underline\pic^\infty_E[\bar g,g]$ to
$\grv\times \grv$ which sends the $S$-valued point $(Y,X,\pi,\bar
x,\bar y,t_{\bar x},t_{\bar y},L,\phi_{\bar y})$ to the following
pair of submodules as a point of $\gr(W)\times \gr(W)$
    {\small
    $$
    \left( \; t_{\bar y}\left(\limil{m}
    (p_* \o_Y(m\cdot \pi^{-1}(\bar x))\right)
    \, , \,
    (t_{\bar y}\circ \phi_{\bar y})\left(\limil{m}
    (p_* L(m\cdot \pi^{-1}(\bar x))\right)\;\right)
    $$}
where $p: Y\times S\to S$ is the projection.

It can be shown that this map is injective and that the image is
contained in the subscheme  $Z\subset \gr(W)\times \gr(W)$ of
those pairs $({\mathcal A},{\mathcal L})$ verifying
    $$
    \o_S\subset {\mathcal A}
    \quad ,\quad
    {\mathcal A}\cdot{\mathcal A}\subseteq {\mathcal A}
    \quad ,\quad
    {\mathcal A}\cdot {\mathcal L}\subseteq {\mathcal L} \ .
    $$

Applying the converse construction of the Krichever correspondence
to the algebra ${\mathcal A}$ we obtain a curve ${\mathcal Y}\to
Z$. Then, consider the subscheme $Z'\subset Z$ defined by the
points $z\in Z$ such that ${\mathcal Y}_z$ is smooth.

Now, we claim that if $({\mathcal A},{\mathcal L})$ belongs to
$Z'$, then ${\mathcal A}$ can be obtained from ${\mathcal L}$.
Indeed, it will be shown that ${\mathcal A}$ is the stabilizer of
${\mathcal L}$.

Consider $({\mathcal A},{\mathcal L})\in Z'$ and let $A_{\mathcal
L}$ denote the stabilizer of $\mathcal L$, that is, the
$\o_S$-algebra
    $$
    A_{\mathcal L}\,:=\,
    \{ w\in \w W_S\text{ such that } w\cdot {\mathcal L}\subseteq {\mathcal
    L}\} \ ,
    $$
Since ${\mathcal A}_z$ corresponds to a smooth curve for all $z\in
Z'$ and ${\mathcal A}\subseteq A_{\mathcal L}$ are points of
$\gr(W)$, then $A_{\mathcal L}$ is a finite ${\mathcal A}$-module
such that ${\mathcal A}_z = (A_{\mathcal L})_z$ for all $z\in Z'$.
Therefore we have that ${\mathcal A}= A_{\mathcal L}$.

Now, one checks that the image of $Z'$ by the projection onto the
second factor,  $\gr(W)\times \gr(W)\to \gr(W)$, represents
$\underline\pic^\infty_E[\bar g,g]$.
\end{pf}

\begin{rem}
For $\chi\in\Z$, the subfunctor of $\underline\pic^\infty_E[\bar
g,g]$ consisting of those points such that $L$ has Euler-Poincar\'{e}
characteristic equal to $\chi$ is representable by the subscheme
$\pic^\infty_E[\bar g,g]\cap \gr^{\chi}(W)$.
\end{rem}

Since $\Gamma_W$ represents the group of invertible elements of
$W$ and $G_V^W$ is a group of automorphisms of $W$ as an algebra,
one has a canonical actions of $G_V^W$ on $\Gamma_W$. Therefore,
it makes sense to consider the semidirect product $G_V^W\ltimes
\Gamma_W$ as follows
    $$
    (g_2,\gamma_2) (g_1,\gamma_1) \,:=\,
    (g_2 g_1,g_1^{-1}(\gamma_2) \gamma_1)
    $$
and the action of $G_V^W\ltimes \Gamma_W$ on the Grassmannian
induced by the action on $W$
    $$
    (g,\gamma) w \,:=\, g(\gamma\cdot w)
    $$

\begin{thm}\label{thm:GGammaactPic}
There are canonical actions of the groups $G_V^W$, $\Gamma_W$ and
$\Gamma_W\ltimes G_V^W$ on $\pic^\infty_E[\bar g,g]$.

Moreover, the action of $G_V^W\ltimes \Gamma_W$ is locally
transitive.
\end{thm}

\begin{pf}
Let
    $$
    \Psi\colon \pic^\infty_E[\bar g,g]\longrightarrow \H^\infty_E[\bar g,g]
    $$
be the forgetful morphism. Let us consider a rational point $p\in
\pic^\infty_E[\bar g,g]$ corresponding to the geometric data
$(Y,X,\pi,\bar x,\bar y,t_{\bar x},t_{\bar y},L,\phi_{\bar y})$.
Let $A:=t_{\bar y}(H^0(Y-\bar y,\o_Y))\in \H^\infty_E[\bar g,g]$
and let $U:= t_{\bar y}(\phi_{\bar y}(H^0(Y-\bar y,L)))\in
\pic^\infty_E[\bar g,g]$.

Considering the map induced b $\Psi$ at the level of tangent
spaces and recalling Theorem~\ref{thm:GonHloctran}, one easily
gets the following diagram
    {\small
    $$
    \xymatrix@C=12pt{
     0 \ar[r] & \lie \Gamma_W \ar[r] \ar@{->>}[d] &
    \lie (\Gamma_W\ltimes G_V^W)
    \ar[d] \ar[r] & \lie G_V^W  \ar@{->>}^{\psi}[d] \ar[r] & 0 \\
    0 \ar[r] &  \hom_{A-\text{mod}}(U,W/U)  \ar[r] & T_p\pic^\infty_E[\bar g,g]  \ar[r] &
    \der(A,W/A)^{\tr} \ar[r] & 0 }
    $$}
The snake lemma implies that the middle vertical arrow is
surjective. We conclude by similar ideas as in the proof of
Theorem~\ref{thm:GonHloctran}.
\end{pf}

\begin{rem}
In a future work and following ideas of \cite{BZF} we aim at
studying how the deformations of a point $\pic^\infty_E[\bar g,g]$
under the action of $G_V^W\ltimes \Gamma_W$ can be interpreted as
isomonodromic deformations.
\end{rem}

%%%%%%%%%%%%%%%%%%%%%%%%%%%%%%%%%%%%%%%%%%


\vskip2truecm
\begin{thebibliography}{MMMi}


\bibitem[BF]{BZF} Ben-Zvi, D.; Frenkel, E., ``Geometric
Realization of the Segal--Sugawara Construction'',
{\verb"math.AG/0301206"}.

\bibitem[ELSV]{elsv} Ekedahl, T.; Lando, S.; Shapiro,
M.; Vainshtein, A., ``Hurwitz numbers and intersections on moduli
spaces of curves'', Invent. Math. {\bf 146} (2001), no. 2, pp.
297--327.

\bibitem[F]{Fulton} Fulton, W.;
``Hurwitz schemes and irreducibility of moduli of algebraic
curves'', Ann. of Math. (2) {\bf 90} (1969), 542--575.

\bibitem[LO]{Levin} Levin, A. M.; Olshanetsky, M. A.,
``Hierarchies of isomonodromic deformations and Hitchin systems'',
Moscow Seminar in Mathematical Physics, 223--262, Amer. Math. Soc.
Transl. Ser. 2, 191, Amer. Math. Soc., Providence, RI, 1999.

\bibitem[MP1]{Hurwitz} Mu\~{n}oz Porras, J. M.; Plaza Mart\'{\i}n, F. J.,
``Equations of Hurwitz schemes in the infinite Grassmannian'',
preprint math.AG/0207091.

\bibitem[MP2]{automorphism} Mu\~noz Porras, J.M.; Plaza Mart\'{\i}n, F.J.,
``Automorphism group of $k(\!(t)\!)$: applications to the bosonic
string'', Commun. Math. Phys. {\bf 216} (2001), pp. 609--634.

\bibitem[OP]{OP}  Okounkov, A.;,  Pandharipande, R.,
``Gromov-Witten theory, Hurwitz theory, and completed
cycles'', math.AG/0204305



\end{thebibliography}
\end{document}